\documentclass[12pt,reqno,a4paper]{amsart}
\usepackage{latexsym,amsmath,amsfonts,amssymb,amsthm}
\theoremstyle{plain}
\newtheorem{Thm}{Theorem}

\theoremstyle{definition}
\newtheorem*{Ack}{Acknowledgment}
\theoremstyle{remark}
\newtheorem*{Rem}{Remark}
\def\Z{\mathbb Z}

\def\Q{\mathbb Q}
\def\jacob #1#2{\genfrac{(}{)}{}{}{#1}{#2}}

\begin{document}

\title{$q$-analogue of Wilson's theorem}
\author{Robin Chapman}
\address{Department of Mathematics, University of Bristol, Bristol,
BS8 1TW, UK} \email{marjc@bris.ac.uk}
\author{Hao Pan}
\address{Department of Mathematics, Nanjing University, Nanjing
210093, People's Republic of China}
\email{haopan79@yahoo.com.cn}
\begin{abstract}We give $q$-analogues of Wilson's theorem for
the primes congruent to $1$ and $3$ modulo $4$ respectively.
Also $q$-analogues of two congruences due to Mordell and Chowla
are established.
\end{abstract} 
\maketitle

\section{Introduction}
\setcounter{equation}{0} \setcounter{Thm}{0} \setcounter{Lem}{0}
\setcounter{Cor}{0}

For arbitrary positive integer $n$, let
$$
[n]_q=\frac{1-q^n}{1-q}=1+q+\cdots+q^{n-1}.
$$
Clearly $\lim_{q\to 1}[n]_q=1$, so we say that $[n]_q$ is a
$q$-analogue of the integer $n$. Supposing that $a\equiv
b\pmod{n}$, we have
$$
[a]_q=\frac{1-q^a}{1-q}
=\frac{1-q^b+q^b(1-q^{a-b})}{1-q}\equiv\frac{1-q^b}{1-q}=[b]_q\pmod{[n]_q}.
$$
Here the above congruence is considered over the polynomial ring $\Z[q]$
in the variable $q$ with integral coefficients. Also $q$-analogues
of some arithmetical congruences have been studied in \cite{Sa,
An, Pa, PS}.

Let $p$ be a prime. The well-known Wilson's theorem states that
$$
(p-1)!\equiv 1\pmod{p}.
$$
Unfortunately, in general,
$$
\prod_{j=1}^{p-1}[j]_q\not\equiv -q^n\pmod{[p]_q}
$$
for any integer $n$. However, we have the following $q$-analogue
of Wilson's theorem for a prime $p\equiv 3\pmod{4}$.
\begin{Thm}
\label{T1} Suppose that $p>3$ is a prime and $p\equiv 3\pmod{4}$.
Then
\begin{equation} \prod_{j=1}^{p-1}[j]_{q^{j}}\equiv
-1\pmod{[p]_q}. \end{equation}
\end{Thm}
In \cite{Mo} (or see \cite[Theorem 8]{UW}), Mordell proved that if
$p>3$ is a prime and $p\equiv 3\pmod{4}$ then
\begin{equation}
\label{mc} \bigg(\frac{p-1}{2}\bigg)!\equiv
(-1)^{\frac{h(-p)+1}{2}}\pmod{p},
\end{equation} where $h(-p)$ is the class number of the quadratic
field $\Q(\sqrt{-p})$. Now we can give a $q$-analogue of
(\ref{mc}).
\begin{Thm}
\label{T2}
 Let $p>3$ be a prime with $p\equiv
3\pmod{4}$. Then
\begin{equation}
\prod_{j=1}^{(p-1)/2}[j]_{q^{16j}}\equiv
(-1)^{\frac{h(-p)+1}{2}}q\pmod{[p]_q}. \end{equation}
\end{Thm}
The case $p\equiv 1\pmod{4}$ is a little complicated. Let
$\jacob{\cdot}{p}$ denote the Legendre symbol modulo $p$. 
By definition,
for any $a$ prime to $p$, $\jacob{a}{p}=1$ or $-1$
according to whether $a$ is a quadratic residue modulo $p$. Let
$\varepsilon_p>1$ and $h(p)$ be the fundamental unit and the class
number of $\Q(\sqrt{p})$ respectively.
\begin{Thm}
\label{T3}
Suppose that $p$ is a prime and $p\equiv
1\pmod{4}$. Then
\begin{equation}
\label{qwp1}
\prod_{j=1}^{p-1}[j]_{q^{j}}\equiv A+B\sum_{j=1,
\jacob{j}{p}=-1}^{p-1}q^j\pmod{[p]_q},
\end{equation}
where
$$
A=\frac{\varepsilon_p^{2h(p)}+\varepsilon_p^{-2h(p)}}{2}
+\frac{\varepsilon_p^{2h(p)}-\varepsilon_p^{-2h(p)}}{2\sqrt{p}}
\text{\quad and\quad
}B=\frac{\varepsilon_p^{2h(p)}-\varepsilon_p^{-2h(p)}}{\sqrt{p}}.
$$
\end{Thm}
Write $\varepsilon_p=(u_p+v_p\sqrt{p})/2$ where $u_p, v_p$ are
positive integers with the same parity. Clearly $u_p^2-pv_p^2=\pm
4$ since $\varepsilon_p$ is an unit. Letting $q\to 1$ in
(\ref{qwp1}), we obtain that
$$
-1\equiv(p-1)!\equiv
A+\frac{B(p-1)}{2}\equiv\frac{\varepsilon_p^{2h(p)}+\varepsilon_p^{-2h(p)}}{2}
\equiv\frac{u_p^{2h(p)}}{2^{2h(p)}}\pmod{p}.
$$
It follows that $h(p)$ is odd and the norm of $\varepsilon_p$ is
always $-1$, i.e., $u_p^2-pv_p^2=-4$.

In \cite{Cho} (or see \cite[Theorem 9]{UW}),
Chowla extended Mordell's result (\ref{mc}) for
$p\equiv 1\pmod{4}$. Let $h(p)$ and
$\varepsilon_p=(u_p+v_p\sqrt{p})/2$ be defined as above. Then
Chowla proved that
\begin{equation}
\label{cc} \bigg(\frac{p-1}{2}\bigg)!\equiv
\frac{(-1)^{\frac{h(p)+1}{2}}u_p}{2}\pmod{p}.
\end{equation}
Now we have the following $q$-analogue of Chowla's congruence:
\begin{Thm}
\label{T4} Suppose that $p$ is a prime and $p\equiv 1\pmod{4}$.
Then
\begin{equation}
\label{qcc} \prod_{j=1}^{(p-1)/2}[j]_{q^{j}}\equiv
-Cq-D\sum_{j=1, \jacob{j}{p}=-1}^{p-1}q^{j+1}\pmod{[p]_q},
\end{equation}
where
$$
C=\frac{\varepsilon_p^{h(p)}-\varepsilon_p^{-h(p)}}{2}
+\frac{\varepsilon_p^{h(p)}+\varepsilon_p^{-h(p)}}{2\sqrt{p}}
\text{\quad and\quad
}D=\frac{\varepsilon_p^{h(p)}+\varepsilon_p^{-h(p)}}{\sqrt{p}}.
$$
\end{Thm}
Let us explain why (\ref{qcc}) implies (\ref{cc}).
Letting $q\to 1$ in (\ref{qcc}), it is derived that
\begin{align*}
\bigg(\frac{p-1}{2}\bigg)!\equiv&-\bigg(\frac{\varepsilon_p^{h(p)}
-\varepsilon_p^{-h(p)}}{2}
+\frac{\varepsilon_p^{h(p)}+\varepsilon_p^{-h(p)}}{2\sqrt{p}}
+\frac{p-1}{2}\cdot \frac{\varepsilon_p^{h(p)}
+\varepsilon_p^{-h(p)}}{\sqrt{p}}\bigg)\\
\equiv&-\frac{((u_p+v_p\sqrt{p})/2)^{h(p)}-((-u_p+v_p\sqrt{p})/2)^{h(p)}}{2}\\
\equiv&-\frac{u_p^{h(p)}}{2^{h(p)}}
=-(u_p^2/4)^{\frac{h(p)-1}{2}}\frac{u_p}{2}\equiv
-\frac{(-1)^{\frac{h(p)-1}{2}}u_p}{2}\pmod{p}.
\end{align*}

The proofs of Theorems \ref{T1}-\ref{T4} will be given in the next sections.

\section{Proofs of Theorems \ref{T1} and \ref{T2}}
\setcounter{equation}{0} \setcounter{Thm}{0} \setcounter{Lem}{0}
\setcounter{Cor}{0}

In this section we assume that $p>3$ is a prime and $p\equiv
3\pmod{4}$. Write
$$
\prod_{j=1}^{p-1}[j]_{q^{j}}=\prod_{j=1}^{p-1}\frac{1-q^{j^2}}{1-q^{j}}
\qquad\text{and}\qquad
\prod_{j=1}^{(p-1)/2}[j]_{q^{16j}}
=\prod_{j=1}^{(p-1)/2}\frac{1-q^{16j^2}}{1-q^{16j}}.
$$
Observe that
$$
[p]_q=\frac{1-q^p}{1-q}=\prod_{j=1}^{p-1}(q-\zeta^j)
$$
where $\zeta=e^{2\pi i/p}$. Also we know that
$\sigma_s:\,\zeta\longmapsto\zeta^s$
is an automorphism over $\Q(\zeta)$ provided
that $p\nmid s$. Hence it suffices to show that
$$
\prod_{j=1}^{p-1}\frac{1-\zeta^{j^2}}{1-\zeta^{j}}=-1\qquad\textrm{and}\qquad
\prod_{j=1}^{(p-1)/2}\frac{1-\zeta^{16j^2}}{1-\zeta^{16j}}
=(-1)^{\frac{h(-p)+1}{2}}\zeta.
$$

Let $Q$ and $N$ denote respectively the sets of quadratic residues
and quadratic non-residues of $p$ in the interval $[1,p-1]$. Then
$$
\prod_{j=1}^{p-1}\frac{1-\zeta^{j^2}}{1-\zeta^{j}}
=\frac{\prod_{j=1}^{\frac{p-1}{2}}(1-\zeta^{j^2})^2}
{\prod_{j=1}^{p-1}(1-\zeta^{j})}=\frac{U^2}{UV}=\frac{U}{V},
$$
where
$$
U=\prod_{k\in Q}(1-\zeta^k),\qquad\text{and}\qquad
V=\prod_{k\in N}(1-\zeta^k).
$$
But since $-1$ is a quadratic non-residue modulo $p$,
$$
V=\prod_{k\in Q}(1-\zeta^{p-k})=\prod_{k\in Q}(1-\zeta^{-k})
=U\prod_{k\in Q}(-\zeta^k).
$$
Now
$$
\sum_{k\in Q}k\equiv\sum_{j=1}^{(p-1)/2}j^2=\frac{p(p^2-1)}{24}
\equiv0\pmod{p}
$$
as $p$ is prime to $6$ and so $p^2\equiv1\pmod{24}$. We conclude
that
$$
U/V=(-1)^{\frac{p-1}{2}}\zeta^{-\sum_{k\in Q}k}=-1
$$
as desired, proving Theorem~\ref{T1}.

Now let us begin to prove
\begin{align}
\label{e21}
\prod_{j=1}^{(p-1)/2}\frac{1-\zeta^{16j^2}}{1-\zeta^{16j}}
=\frac{\prod_{j=1}^{(p-1)/2}(1-\zeta^{16j^2})}
{\prod_{j=1}^{(p-1)/2}(1-\zeta^{16j})}=(-1)^{\frac{h(-p)+1}{2}}\zeta.
\end{align}
Clearly the numerator of the left side of (\ref{e21}) is $U$. Let
$$
W=\prod_{j=1}^{(p-1)/2}(1-\zeta^{16j})
$$
denote its denominator. Let $M=\{1,2,\ldots,(p-1)/2\}$. Then
$W=W_+W_-$ where
$$
W_+=\prod_{j\in M\cap Q}(1-\zeta^{16j})\qquad\textrm{and}\qquad
W_-=\prod_{j\in M\cap N}(1-\zeta^{16j}).
$$
Now
$$
W_-=\sum_{M'\cap Q}(1-\zeta^{-16k})=\frac{U}{W_+} \prod_{k\in
M'\cap Q}(-\zeta^{-16k})
$$
where $M'=\{(p+1)/2,\ldots,p-1\}$. We know (see 
\cite[Chapter 5, Section 4, Theorem 3]{BS}) that
\begin{align*}
h(-p)=&\frac{1}{2-\jacob{2}{p}}\sum_{k=1}^{(p-1)/2}\jacob{k}{p}
=\frac{1}{2-\jacob{2}{p}}
\bigg(\frac{p-1}{2}-2\sum_{k\in M\cap N}1\bigg)\\
\equiv&-1-2|M\cap N|=-1-2|M'\cap Q|\pmod{4}.
\end{align*}
Also, we have
\begin{align*}
\frac{p^2-1}{8}=&\sum_{k=1}^{(p-1)/2}k=\sum_{k\in M\cap
Q}k+\sum_{k\in M\cap N}k=\sum_{k\in M\cap Q}k+\sum_{k\in M'\cap Q}(p-k)\\
\equiv&\sum_{k\in Q}k-2\sum_{k\in M'\cap Q}k\equiv-2\sum_{k\in
M'\cap Q}k\pmod{p},
\end{align*}
whence $\sum_{k\in M'\cap Q}16k\equiv 1\pmod{p}$. Thus
$$
\frac{U}{W}=\frac{U}{W_+W_-}=(-1)^{|M'\cap Q|}\prod_{k\in M'\cap
Q} \zeta^{16k}=(-1)^{(1+h(-p))/2}\zeta,
$$
which confirms (\ref{e21}).\qed

\section{Proofs of Theorems \ref{T3} and \ref{T4}}
\setcounter{equation}{0} \setcounter{Thm}{0} \setcounter{Lem}{0}
\setcounter{Cor}{0}

Below suppose that $p$ is a prime congruent to $1$ modulo $4$
and $\zeta=e^{2\pi i/p}$. Let $Q,\ N\subseteq[1,p-1]$
be the sets of quadratic residues
and quadratic non-residues of $p$ respectively. Let
$$
U=\prod_{k\in Q}(1-\zeta^k)\qquad\text{and}\qquad
V=\prod_{k\in N}(1-\zeta^k).
$$
In order to prove Theorem \ref{T3}, we only need to show prove that
\begin{equation}
\label{e31}
A+B\sum_{j\in
N}\zeta^j=\prod_{j=1}^{p-1}\frac{1-\zeta^{j^2}}{1-\zeta^j}=\frac{\prod_{j\in
Q}(1-\zeta^{j})^2}{\prod_{j\in
Q}(1-\zeta^{j})\prod_{j\in
N}(1-\zeta^{j})}=\frac{U}{V}.
\end{equation}
By the analytic class number formula \cite[Chapter 1, Section 4,
Theorem 2]{BS}
$$U=\varepsilon_p^{-h(p)}\sqrt{p}\qquad\textrm{and}\qquad
V=\varepsilon_p^{h(p)}\sqrt{p}.$$ Thus
$U/V=\varepsilon_p^{-2h(p)}=a-b\sqrt{p}$ where
$$2a=\varepsilon_p^{2h(p)}+\varepsilon_p^{-2h(p)}\in\Z
\quad\textrm{and}\quad 2b
=(\varepsilon_p^{2h(p)}-\varepsilon_p^{-2h(p)})/\sqrt{p}\in\Z.$$
Also, by Gauss's formula for the quadratic Gauss sum
$$
\sqrt{p}=\sum_{j=1}^{p-1}\jacob{j}{p}\zeta^j
=\sum_{j=1}^{p-1}\zeta^j-2\sum_{j\in N}\zeta^j=-1-2\sum_{j\in N}\zeta^j.
$$
Hence
$$
\frac{U}{V}=a+b\bigg(1+2\sum_{j\in N}\zeta^j\bigg),
$$
which is clearly equivalent to (\ref{e31}).

\begin{Rem}
in \cite{Cha}
the first author used products like $\prod_{j\in N}(1-\zeta^j)$
to study determinants built from Legendre symbols.
\end{Rem}

Let us now consider the product
$$
\prod_{j=1}^{(p-1)/2}\frac{1-\zeta^{16j^2}}{1-\zeta^{16j}}=\frac{U}{\Pi_{16}}
$$
where
$$\Pi_r=\prod_{j=1}^{(p-1)/2}(1-\zeta^{rj}).$$
When $r$ and $s$ are prime to $p$, we have $\Pi_{rs}=\sigma_s(\Pi_r)$
where $\sigma_s$ is the automorphism of $\Q(\zeta)$ mapping $\zeta$
to $\zeta^s$. It turns out to be convenient to compute $\Pi_{16}$
as $\sigma_4(\Pi_4)$. As $p\equiv1\pmod{4}$ we know that
$U=\varepsilon_p^{-h(p)}\sqrt{p}$. For each $r$ prime to~$p$,
$$
|\Pi_r|^2=\Pi_r\Pi_{-r}=\prod_{j=1}^{p-1}(1-\zeta^j)=p,
$$
so $|\Pi_r|=\sqrt{p}$.
Now
\begin{align*}
\frac{\Pi_4}{|\Pi_4|}
&=\prod_{j=1}^{(p-1)/2}\frac{1-\zeta^{4j}}{|1-\zeta^{4j}|}
=\prod_{j=1}^{(p-1)/2}(-\zeta^{2j})\frac{2i\sin(4\pi j/p)}
{|2i\sin(4\pi j/p)|}\\
&=(-1)^M\prod_{j=1}^{(p-1)/2}(-i\zeta^{2j})
=(-1)^M(-i)^{(p-1)/2}\zeta^{(p^2-1)/4}\\
&=(-1)^{(p-1)/4+M}\zeta^{(p^2-1)/4}
\end{align*}
where
$$
M=|\{1\leqslant j\leqslant(p-1)/2:\,\sin(4\pi j/p)<0\}|.
$$
Note that when $0<j<p/2$, $\sin(4\pi j/p)<0$ if and only if $p/4<j<p/2$. So
$M=(p-1)/4$, and
$$
\Pi_4=\zeta^{(p^2-1)/4}\sqrt{p}.
$$
Also
$$
\sigma_4(\sqrt{p})
=\sigma_4\bigg(\sum_{j=1}^{p-1}\jacob{j}{p}\zeta^j\bigg)
=\sum_{j=1}^{p-1}\jacob{j}{p}\zeta^{4j}
=\sum_{j=1}^{p-1}\jacob{j}{p}\zeta^{j}=\sqrt{p}.
$$
Thus
$$
\Pi_{16}
=\sigma_4(\zeta^{(p^2-1)/4}\sqrt{p})=\zeta^{p^2-1}\sqrt{p}=\zeta^{-1}\sqrt{p}.
$$

Assume that $\varepsilon_p^{h(p)}=(c+d\sqrt{p})/2$ where $c$ and $d$
are integers with the same parity. Recall that
the norm of $\varepsilon_p$ is $-1$ and $h(p)$ is odd.
Hence $\varepsilon_p^{-h(p)}=(-c+d\sqrt{p})/2$.
So
$$
c=\varepsilon_p^{h(p)}-\varepsilon_p^{-h(p)}\qquad\textrm{and}\qquad
d\sqrt{p}=\varepsilon_p^{h(p)}+\varepsilon_p^{-h(p)}.
$$
As $\sqrt{p}=-1-2\sum_{j\in N}\zeta^j$, we have
$$
\varepsilon_p^{-h(p)}=\frac{-c+d\sqrt{p}}{2}
=-\frac{c+d}{2}-d\sum_{j\in N}\zeta^{j}.
$$
Therefore
$$
\prod_{j=1}^{(p-1)/2}\frac{1-\zeta^{16j^2}}{1-\zeta^{16j}}
=\frac{U}{\Pi_{16}}=\frac{\varepsilon_p^{-h(p)}\sqrt{p}}{\zeta^{-1}\sqrt{p}}
=-\frac{(c+d)\zeta}{2}-d\sum_{j\in N}\zeta^{j+1}.
$$
This implies that
$$
\prod_{j=1}^{(p-1)/2}[j]_{q^{16j}}\equiv
-\frac{c+d}{2}q-d\sum_{j\in N}q^{j+1}\pmod{[p]_q}.
$$
We are done.\qed

\begin{Ack} The second author thanks Professor Zhi-Wei Sun for his
helpful suggestions on this paper. \end{Ack}

\end{document}